\pageno=-1
\magnification\magstep1

\newif\ifproofmode
%\proofmodetrue

\font\bigfont=cmr17
\font\small=cmsl10 at 10 true pt
\font\teneuf=eufm10
\font\seveneuf=eufb7
\font\fiveeuf=eufb5
\font\msbm=msbm10
\font\sevenmsbm=msbm7
\font\fivemsbm=msbm5
\font\msam=msam10
\font\sevenmsam=msam7
\font\fivemsam=msam5

\def\circle#1{}
\ifproofmode 
\else
\font\circle=lcircle10
\fi

\catcode`\@11

\advance\topskip \baselineskip
\abovedisplayskip0.5\abovedisplayskip
\belowdisplayskip0.5\belowdisplayskip

\hfuzz5truept
\let\BEGINsection\beginsection % why? because ispell thinks thatthis
			       % is f***ing Latex, in which \begin...
			       % always needs {} and and end ...

\def\,{\ifmmode\mskip\thinmuskip\else$\,$\fi}

%\footline{\tt \jobname\hfil\tenrm\folio
%\hfil\number\day.\number\month.\number\year}
%% a few essential macros

\def\proof{\par\noindent Proof}
\def\proc#1:{\advance\thmcount1 
\smallskip\noindent{\bf \number\thmcount. #1:}}
\def\<#1>{\hbox{$\langle #1\rangle$}}
\newcount\thmcount

\def\cut{\cap}

\def\Limpl{\Rightarrow}
\def\liff{\leftrightarrow}
\def\Liff{\Leftrightarrow}
\def\cf{{\it cf}}

\def\BEGINdent{\par\begingroup\parskip0cm\par
	\advance\parindent1cm\advance\rightskip0.8cm}
\def\ENDent{\par\endgroup\par}
\def\ite #1 {\item{(#1)}}
 \def\forces{\mathrel{\mathrel\Vert\joinrel\mathrel-}}

\def\Limpl{\Rightarrow}
\def\Levy{{\rm Levy}}

\newcount\skewfactor
\def\mathunderaccent#1#2{\let\theaccent#1\skewfactor#2
\mathpalette\putaccentunder}
\def\putaccentunder#1#2{\oalign{$#1#2$\crcr\hidewidth
\vbox to.2ex{\hbox{$#1\skew\skewfactor\theaccent{}$}\vss}\hidewidth}}
\def\skewname#1#2{\mathunderaccent\tilde{#1}{#2}}
\def\assign#1 #2 {
\edef\next{\noexpand\def
\csname name#1\endcsname{\noexpand\skewname{#2}{#1}}}%
\next}
\def\name#1{\skewname{-3}{#1}}

\assign A -1
\assign G -3 
\assign P -8 
\assign Q -2 
\assign S -3 
\assign T -6 
\assign a -1
\assign b -2 
\assign c 0 
\assign d -2 
\assign e 0 
\assign f -6
\assign g -4 
\assign h -2 
\assign i -1 
\assign l -2
\assign p -6
\assign q 0
\assign t -1 
\assign x -1 
\assign y -3 
\assign z -1

\newfam\frakturfam
\textfont\frakturfam\teneuf
\scriptfont\frakturfam\seveneuf
\scriptscriptfont\frakturfam\fiveeuf

\def\frak{\fam\frakturfam\teneuf}

\newfam\msyfam
\textfont\msyfam=\msbm
\scriptfont\msyfam=\sevenmsbm
\scriptscriptfont\msyfam=\fivemsbm

\newfam\msxfam
\textfont\msxfam=\msam
\scriptfont\msxfam=\sevenmsam
\scriptscriptfont\msxfam=\fivemsam

\def\Bbb{\relax\ifmmode\fam\msyfam\else\message{Bbb not allowed in
text}\fi}

\def\hexnumberat#1{\ifcase#1 0\or 1\or 2\or 3\or 4\or 5\or 6\or 7\or 8\or
 9\or A\or B\or C\or D\or E\or F\fi}
\edef\msxm{\hexnumberat\msxfam}
\edef\msym{\hexnumberat\msyfam}

\mathchardef\upharpoonright"3\msxm16

\def\on{{\upharpoonright}}

\mathchardef\beth"0\msym69

\newcount\skewfactor
\def\mathunderaccent#1#2 {\let\theaccent#1\skewfactor#2
\mathpalette\putaccentunder}
\def\putaccentunder#1#2{\oalign{$#1#2$\crcr\hidewidth
\vbox to.2ex{\hbox{$#1\skew\skewfactor\theaccent{}$}\vss}\hidewidth}}
\def\name{\mathunderaccent\widetilde-1 }

%% The tilde is not quite in the exact right place.  I will think of
%% something. 

\def\qedf#1{}
\def\qed#1{\hfill\rlap{\hskip\rightskip$\smiley_{#1}$}$\phantom{\smiley{#1}}$%
}
\def\QED#1{\hfill\rlap{\hskip\rightskip$\SMILEY_{#1}$}$\phantom{\smiley{#1}}$%
}

\setbox0=\hbox{~~~~~}
\setbox1=\hbox to \wd0{\hfill$\scriptstyle\smile$\hfill} % mouth
\setbox2=\hbox to \wd0{\hfill$\cdot\,\cdot$\hfill} % eye 
\setbox3=\hbox to \wd0{\hfill\hskip4.8pt\circle i\hskip-4.8pt\hfill} % head
%\setbox4=\hbox to \wd0{\hfil\vrule width 0.2pt 
%                       height 3.5pt depth -1pt\hfil}% nose
\setbox5=\hbox to  \wd0{\hfill--\hfill}

\wd1=0cm
\wd2=0cm
\wd3=0cm
\wd4=0cm
\wd5=0cm

\newbox\SMILEBOX
\setbox\SMILEBOX \hbox {\lower 0.4ex\copy1
                 \raise 0.3ex\copy2
                 \raise 0.5ex\copy3
                \copy4
                \copy0{}}
\def\SMILEY{\leavevmode\copy\SMILEBOX}

\newbox\smilebox
\setbox\smilebox \hbox {\lower 0.4ex\box5
                 \raise 0.3ex\box2
                 \raise 0.5ex\box3
                \box4
                \box0{}}
\def\smiley{\leavevmode\copy\smilebox}

\newlinechar`\^^J
\catcode`@=11
\newwrite\@unused
\newwrite\mgfile

\def\@ifstar#1#2{\@ifnextchar *{\def\@tempa*{#1}\@tempa}{#2}}

\def\@ifnextchar#1#2#3{\let\@tempe #1\def\@tempa{#2}\def\@tempb{#3}\futurelet
    \@tempc\@ifnch}
\def\@ifnch{\ifx \@tempc \@sptoken \let\@tempd\@xifnch
      \else \ifx \@tempc \@tempe\let\@tempd\@tempa\else\let\@tempd\@tempb\fi
      \fi \@tempd}

\newcount\secno\newcount\theono
 
\def\BAremark{0.3} % (2)
\def\firstCBA{0.4} % (2)
 % (4)
\def\PIisFFA{1.4} % (5)
\def\commuting{1.8} % (7)
\def\fhasnofixpoints{1.10} % (7)
\def\lastfact{1.11} % (8)
\def\twotrees{1.13} % (9)
\def\regrem{2.3} % (12)
\def\inacc{2.4} % (13)
 % (13)
\def\Hchi{2.7} % (13)
\def\properisabsolute{2.8} % (14)
\def\findsmallproper{2.9} % (14)
\def\smallforcing{2.10} % (14)
\def\whatisPT{3.2} % (17)
\def\mitchelfact{3.3} % (17)
 % (18)
\def\equicon{4.1} % (19)
\def\absolute{4.2} % (19)
\def\howtodefineit{4.4} % (21)
\def\freeofro{4.5} % (21)
\def\secondCBA{4.6} % (21)
\def\whyWF{4.7} % (22)
\def\coveringlemma{4.8} % (23)
\def\WFandTheory{4.11} % (25)
\def\thinksphi{4.13} % (26)
\def\isacardinal{4.14} % (26)
\def\nuisthere{4.15} % (26)
\def \PFApplications{Ba1}
\def \IteratedForcing{Ba2}
\def \fuchino{Fu}
\def \Mitchell{Mi}
\def \EhrenfeuchtConjecture{Sh 56}
\def \LXXIII{Sh\penalty \@M \ 73}
\def \ProperForcing{Sh\penalty \@M \ b}
\def \ProperImproper{Sh\penalty \@M \ f}
\def \todo{To}

\openin\mgfile \jobname.mg
\ifeof\mgfile \message{No file \jobname.mg}
	\else\closein\mgfile\relax\input \jobname.mg\fi
\relax
\openout\mgfile=\jobname.mg

\def\@nofirst#1{}

\def\neusection{\advance\secno by 1\relax \theono=0\relax}
\def\neuchap{\secno=0\relax\theono=0\relax}

\neuchap

\def\labelitshort{\global\advance\theono by 1 
\global\edef\thistheorem{\number\secno.\number\theono}}

\def\labelit#1{\global\advance\theono by 1%
             \global\edef#1{%
             \number\secno.\number\theono}%
\let\thistheorem#1%
             \write\mgfile{\@definition{#1}}%
}

\def\proc{\@ifstar{\npro}{\spro}}
\def\@wemustnotuse{thistheorem}

\def\spro#1:{%
\labelitshort%
\smallbreak\noindent%
\@markit{\thistheorem}{}%
{\bf\ignorespaces#1:}}

\def\npro#1 #2:{%
\def\temp@{#1}\ifx\temp@\@wemustnotuse
\immediate\write\@unused{^^J^^JWARNING: Please do not use the label 
``thistheorem''^^J^^J}%
\fi
\expandafter\labelit\csname#1\endcsname%
\smallbreak\noindent%
\expandafter\@markit\csname#1\endcsname{#1}%
{\bf \ignorespaces #2:}}

\def\@definition#1{\string\def\string#1{#1}
\expandafter\@nofirst\string\%
(\the\pageno)}

\def\@markit#1#2{% should be called only at beginning of paragraph
\ifproofmode\llap{{\tt #2\ }}\fi%
{\bf #1\ }%
}
 
\def\labelcomment#1{\write\mgfile{\expandafter
		\@nofirst\string\%---#1}} 

\catcode`@=12

%%%%%%%% file labelit.tex end %%%%%%%%

%\proofmodefalse

\def\setup#1#2{\def\wherereflect{#1}\def\really{#2}%
\futurelet\klammer\investigateit}
\def\investigateit{\ifx\klammer[%]
  \def\next{\waitforclose}\else\def\next{\really}\fi\next}
\def\waitforclose[#1]{\edef\wherereflect{\wherereflect[#1]}\really}

\openup 2pt
\def \barchi {{\chi'}}
\def\logand{\,\mathbin{\&}\,}
\def\B{{\frak B}}
\def\dom{{\rm dom}}\def\range{{\rm rng}}

\def\VEC#1{#1_1, \ldots, #1_k}

\def\txt#1{\def\inhbox{#1}\futurelet\maybeaquote\checkifquote}
\def\checkifquote{\ifx\reallyaquote\maybeaquote\let\next\gobbleandmore%
\else\def\next{\hbox{\inhbox}}\relax\fi\next}
\def\gobbleandmore#1{\edef\inhbox{\inhbox#1}%
\futurelet\maybeaquote\checkifquote}
\let\reallyaquote'
	
\def\\{\hfil\break}
\def\E{{\cal E}}
\def\A{{\frak A}}
\def\B{{\frak B}}
\def\C{{\frak C}}
\def\M{{\frak M}}
\def\N{{\frak N}}
\def\L{{\cal L}}
\def\supp{{\rm supp}}
\def\BFA{{\rm BFA}}
\def\BPFA{{\rm BPFA}}
\def\iso{\simeq}
\def\Suc{{\rm Suc}}
\def\concat{{}^\frown}
\def\rk{{\rm rk}}

\baselineskip1.6\baselineskip
\parskip=0.2\baselineskip
\advance\parskip0pt plus 1pt
\parindent0cm

\medskip

{{

%\vglue 3cm

\centerline{\bigfont
The Bounded Proper Forcing Axiom
}}
\baselineskip0.8\baselineskip
\bigskip
\centerline{Revised version, January 1994}
\bigskip
\bigskip

\centerline{\bf Martin
GOLDSTERN\rm\footnote{$^1$}{\openup-7pt
The authors thank the  DFG (grant  Ko 490/7-1)
and  the Edmund Landau Center for research in Mathematical
        Analysis, supported by the Minerva Foundation (Germany)%
\parskip0pt\endgraf\vskip-2\baselineskip
}
}

\centerline{Free University of Berlin}

\bigskip

\centerline{\bf
Saharon SHELAH\rm$^{1,}$\footnote{$^2$}{Publication 507}}

\centerline{Hebrew University of Jerusalem}

\bigskip
\bigskip

{\leftskip0.15\hsize\rightskip0.15\hsize\noindent ABSTRACT. 
\rm \parskip0pt
The bounded proper forcing axiom BPFA is the statement that 
 for any family
of $\aleph_1$ many maximal antichains of a proper forcing notion, each
of size $\aleph_1$, there is a directed set meeting all these
antichains. 

A regular cardinal $\kappa$ is called ${\Sigma}_1$-reflecting, if for any
regular cardinal $\chi$, for all formulas $\varphi$, 
``$H(\chi) \models `\varphi$'\,'' implies ``$\exists {\delta} {<}
\kappa$, $H({\delta}) \models `\varphi$'\,'' 

We show that BPFA is equivalent to the statement that two
nonisomorphic models of size $\aleph_1 $ cannot be made isomorphic by
a proper forcing notion, and we show that the consistency strength of
the bounded proper forcing axiom is exactly the existence of a
${\Sigma}_1$-reflecting cardinal (which is less than the existence of
a Mahlo cardinal). 

We also show that the question of the existence of isomorphisms
between two structures can be reduced to the question of rigidity of a
structure. 

}

\nopagenumbers
\bigskip\bigskip

\eject

}
\headline{\small\hfill Goldstern--Shelah: Forcing Axioms with Small
 Antichains \hfill}

\pageno=1

\BEGINsection{Introduction}

The proper forcing axiom has been successfully employed to decide many
questions in set-theoretic topology and infinite combinatorics.  See
[\PFApplications] for some applications, and [\ProperForcing] 	and
[\ProperImproper] for variants.

In the recent paper [\fuchino], Fuchino investigated the following two
consequences of the proper forcing axiom: 
\BEGINdent
\ite a If a structure $\A$ of size $\aleph_1$ cannot be embedded into
a structure $\B$, then such an embedding cannot be produced by a
proper forcing notion. 
\ite b If two structures $\A$ and $\B$ are not isomorphic, then they
cannot be made isomorphic by a proper forcing notion. 
\ENDent
He  showed that (a) is in fact equivalent to the proper forcing axiom,
and asked if the same is  true for (b).  

In this paper we find a natural weakening of the proper forcing axiom,
the ``bounded'' proper forcing axiom  and show that it is equivalent
to property (b) above.   

We then investigate the consistency strength of this new axiom.  While
the exact consistency strength of the proper forcing axiom is still
unknown (but large, see [\todo]), it turns out that the bounded proper forcing
axiom is equiconsistent to a rather small large cardinal.  

%   The purpose of this paper is to present a variant 
% of PFA which on the one hand is equiconsistent to a
% (mildly) large cardinal, and on the other hand has the full power of 
%  PFA as far as subsets of $\omega _1$ are concerned. 

\bigskip

For notational simplicity we will, for the moment,  only consider
forcing notions which are complete Boolean algebras.    See \firstCBA\  and
\secondCBA.

We begin by recalling the forcing axiom in its usual form:
For a forcing notion $P$, 
FA$(P,\kappa)$ is the following statement: 
\BEGINdent
\item {} Whenever $\< A_i:i< \kappa > $ is a family of maximal
antichains of $P$, then there is a filter $G^* \subseteq P$
meeting all $A_i$. 
\ENDent

If $\name f $ is a $P$-name for a function from $\kappa$ to the
ordinals, we will say that $G^* \subseteq P$ decides $\name f$ if for each
$i<\kappa $ there is a condition $p\in G^*$ and an ordinal $\alpha_i$
such that $p \forces \name f(i) = \alpha_i$. (If $G^*$ is directed,
then this ordinal must be unique, and we will write $\name f[G^*]$ for
the function $i  \mapsto \alpha_i$.)   Now it is easy to see
that the FA$(P,\kappa)$ is equivalent to the following statement: 
\BEGINdent
\item {} Whenever $\name f$ is a $P$-name for a function from $\kappa$
to the ordinals, then there is a filter $G^* \subseteq P$ which
decides $\name f$. 
\ENDent
      
This characterization suggests the following weakening of the forcing
axiom: 

\proc Definition:  Let $P$ be a forcing notion, and let $\kappa$
and ${\lambda}$ be infinite cardinals.
\BEGINdent
\item {} $ \BFA(P,\kappa,{\lambda})$ is the following statement:
	Whenever $\name f$ is a $P$-name for a function from $\kappa$
	to ${\lambda}$ then there is a filter $G^* \subseteq P$
	which decides $\name f$, or equivalently:
\item{} 	Whenever $\<A_i:i<\kappa >$ is a family of maximal antichains
	of $P$, each of size $\le {\lambda}$,  then there is a
	filter $G \subseteq P$ which meets all $A_i$. 

\ENDent

\proc Notation: 
\BEGINdent
\ite 1 $\BFA(P,\lambda) $ is $\BFA(P,\lambda,\lambda)$, and $\BFA(P) $ is
	$\BFA(P,\omega_1)$. 
\ite 2 If $\E$ is a class or property  of forcing notions, we write
	$\BFA(\E)$ for 
	$\forall P{\in} \E\,\, \BFA(P)$, etc. 
\ite 3 $\BPFA$ = the bounded proper forcing axiom =  $\BFA({\rm
	proper})$. 
\ENDent
      
\proc *BAremark Remark: 
 For the class of ccc forcing notions we get nothing new:
$\BFA({\rm ccc}, \lambda)$ is equivalent to  Martin's axiom
MA$(\lambda)$, i.e., FA$({\rm ccc}, \lambda)$. \qed\BAremark

\proc *firstCBA Remark:  If the forcing notion $P$ is not a complete
Boolean 
algebra but an arbitrary poset, then it is possible that $P$ does not
have any small antichains, so it could satisfy the second version of
BFA$(P)$ vacuously.  The problem with the first definition, when
applied to an arbitrary poset, is that a filter on ro($P$) which
interprets the $P$-name (=ro$(P)$-name) $\name f$ does not necessarily
generate a filter on $P$.   So for the moment our official definition
of BFA$(P)$ for arbitrary posets $P$ will be
$$ \BFA(P) \ :\iff \ \BFA({\rm ro}(P))$$
In \howtodefineit\ and \freeofro\  
we will find a equivalent (and more natural?)\ definition $\BFA'(P)$ which
does not explicitly refer to ro$(P)$.

\bigskip\bigskip\bigskip
Contents of the paper:  in section 1 we show that the ``bounded forcing
axiom'' for any  forcing notion $P$ is equivalent
to Fuchino's ``potential isomorphism'' axiom for $P$.
  In section 2 we define
the concept of a ${\Sigma}_1$-reflecting cardinal, and we show that from
a model with such a cardinal we can produce a model for the bounded
proper forcing axiom. In section 3 we describe a (known) forcing notion
which we will use in section 4, where we complement our consistency  result
by showing that a ${\Sigma}_1$-reflecting  cardinal is necessary:  If
BPFA holds, then $\aleph_2$ must be ${\Sigma}_1$-reflecting in $L$. 

\bigskip

Notation: 
We use \SMILEY{} to denote the end of a proof, and we write \smiley{} when
we leave a proof to the reader. 

We will use gothic letters $\A$, $\B$, $\M$, \dots\  for structures
(=models of a first order language),
and the corresponding latin letters $A$, $B$, $M$, \dots \ 
 for the underlying universes.   Thus, a model $\A$ will have the
universe $A$, and if $A' \subseteq A$ then we let $\A'$ be the submodel
(possibly with partial functions) with universe $A'$, etc. 

\neusection
\BEGINsection 1. Fuchino's problem and other applications

Let $\E$ be a class  of forcing notions.

\proc Definition: 
 Let $\A$ and $\B$ be two
structures for the same first order language, and let $\E$ be a class
(or property) of  forcing notions.   We say that $\A$ and $\B$
are $\E$-potentially isomorphic ($\A\iso_\E B$) iff there is a forcing
$P\in \E$ such that 
 $\forces _P $``$\A \iso \B$.''  $\A\iso_P\B$ means
$\A\iso_{\{P\}}\B$. 

\proc Definition:  We say that a structure $\A$ is nonrigid, if it
admits a nontrivial automorphism.  We say that $\A$ is
$\E$-potentially nonrigid, if there is a forcing notion $P\in \E$, 
$\forces_P$ ``$\A$
is nonrigid''.

\proc *PIandPA Definition:
\BEGINdent
\ite 1  PI$(\E, \lambda) $  is the statement: Any two $\E$-potentially
isomorphic structures of size ${\lambda}$ are isomorphic.
\ite 2 PA$(\E, {\lambda})$ is the statement: Any $\E$-potentially nonrigid
structure of size ${\lambda}$
is nonrigid.
\ENDent
      
PI$(\E,{\lambda})$ was defined by Fuchino [\fuchino].  It is clear that
$${\rm FA}(\E, \lambda)\ 
 \Longrightarrow \ 
\BFA(\E, \lambda ) 
\ \Longrightarrow \  
{\rm PI}(\E,{\lambda})
\&
{\rm PA}(\E,{\lambda})
$$ for all $\E$,  and Fuchino asked if PI$(\E,{\lambda})$ implies
FA$(\E, \lambda)$, in particular for the cases 
$\E$=ccc, $\E$=proper and $\E$=stationary-preserving.  

We will show in this section that the the three statments 
BFA, PA and PI are in fact
equivalent.  Hence
in particular PI$({\rm ccc},\lambda)$ is equivalent to MA$(\lambda)$.  

In the next sections we will show that for $\E$=proper, the first
implication cannot 
be reversed, by computing the exact consistency strength of BPFA and
comparing it to the known lower bounds for the consistency strength of
PFA.

\proc *PIisFFA Theorem: For any forcing notion $P$ and for any
${\lambda}$, the following are equivalent: 
\BEGINdent
\item{} 
       PI$(P,{\lambda})$
\item {} PA$(P,{\lambda})$ 
\item {} $\BFA(P, \lambda)$ 
\ENDent

% restriction  closed class $\E$. (We say that $\E$ is restriction-closed
% if $p\in Q\in \E$ implies $\{q\in Q:q \hbox{ is stronger than } p\}\in
% \E$.)
% 
% 
% 
% 
% 
% 
% \proc  Corollary: 
% \BEGINdent
% \ite 1 
%  PI$({\rm ccc},{\lambda}) \ \Leftrightarrow\ MA({\lambda}) $ 
% \ite 2  PI$({\rm proper}, \lambda) \
% \Leftrightarrow \ \BPFA(\lambda)$. \qed\thistheorem
% \ENDent

\proof\ of PI $\Limpl$ PA:   This follows from theorem \twotrees\
below.  Here we will give a shorter  proof under the additional
assumption that we have not only PI($P$) but also PI($P_p$) for all
$p\in P$, where $P_p$ is the set of all elements of $P$ which are
stronger that $p$:  

  Let $\M$ be a potentially nonrigid
structure.    So there is a a $P$-name
$\name f$ such that 
$$\forces_P \hbox{``$\name f$ is a nontrivial automorphism of $\M$''}$$
We can find a condition $p\in P$ and two elements $a\not=b$ of $\M$
such that 
$$p\forces_P \hbox{``$\name f(a)=b$''}$$
Since we can replace $P$ by $P_p$, we may assume that
$p$ is the weakest condition of $P$.  So we have that $(\M,a)$ and
$(\M,b)$ are potentially isomorphic.  Any isomorphism from $(\M,a)$ to
$(\M,b)$ is an automorphism of $\M$ mapping $a$ to $b$, so we are
done. 
\QED{PI \Limpl PA}

% 
% \proc  Fact: 
% \BEGINdent
% \ite 1 $\supp(w)$ is a finite subset of $X$. 
% \ite 2 Let  $g:F(X) \to F(Y) $ be a homomorphism, $w\in F(X)$.  If
% $\supp(g(w))\not=\emptyset$ then $\supp(w)\not=\emptyset$. 
% \ENDent
      
% \proc *denseequivalent Remark:  If $Q$ is a dense subset of $P$, then
% PA$(P,\lambda)$ is equivalent to  PA$(Q,\lambda)$, and  
% $\BFA(P,\lambda)$ is equivalent to $\BFA(Q,\lambda)$. 
% 
% Proof:  $P$ and $Q$ generate the same generic extension.  In
% particular there is a canonical correspondence between $P$-names and
% $Q$-names. \QED\denseequivalent
% 
\bigskip

We will now describe the framework of the proof 
of the second part of our
theorem: PA $\Limpl  $ BFA.
 We start with a forcing notion $P$.   Recall that (for the moment)
all our forcing notions are a complete Boolean
algebras.
 Fix a small family of small
antichains. 
 Our structure will consist of a disjoint union of the free
groups generated by the antichains.  On each of the free groups the
translation by an element of the corresponding antichain will be a
nontrivial automorphism, and if all these translating elements are
selected from the antichains by a directed set, then the union of
these automorphisms will be an automorphism of the whole structure.
We will also ensure that ``essentially'' these are the only
automorphisms, so every automorphism will define a sufficiently
 generic set.

\proc Definition: For any set $X$ let $F(X)$ be the free group on the
generators $X$, and for $w\in F(X)$ define 
$\supp(w) = \bigcap \{ Y \subseteq X : w \in \langle Y \rangle\}$,
i.e.\ $\supp(w)$ is the set of elements $x$ of $X$ which occur (as $x$
or as $x^ {-1}$) in the reduced representation of $w$. 
(If you prefer, you can change the proof below by using the free
abelian group generated by $X$ instead of the free group,  or the free
abelian group of order 2,  \dots).

\proc Setup: 
Let $P$ be a complete Boolean algebra, and 
let $(A_i:i \in I )$ be a
system of ${\lambda} $ many 
maximal antichains of size $\lambda$.  We may assume that this is a directed
system, i.e., for any $i,j \in I $ there is a $k\in I$ such that
$A_k$ refines both $A_i$ and $A_j$.  So if we  write $i<j$ for ``$A_j$
refines $A_i$'', $(I,<)$ becomes a partially ordered upwards directed
set.  
(We say that $A$ refines $B$ if each element of $A$ is stronger than
some unique element of $B$, or in the Boolean sense if there is a
partition
  $ A = \bigcup\limits_{b\in B }^{\hskip-.5pt\vrule width 1pt height1pt} A_b$ of
the set $A$  
satisfying $\forall b\in B$ $\sum_{a\in A_b} a = b$.) 
 
Assuming PA$(P,{\lambda})$, we will find a filter(base) meeting all
the sets $A_i$. 

\proc Definition: 
\BEGINdent
\ite  a 
Let $(F_i, *)$  be the free group generated by $A_i$, and let $M$ be the
disjoint union of the sets $F_i$.   
\ite b For $i\in I$, $z\in F_i$ let 
$$ R_{i,z} = \{ (y, z * y): y \in F_i\}$$
\ite c If $i < j$, then there is a ``projection'' function $h^j_i$ 
from $A_j$ to $A_i$: For $p\in A_j$, $h^j_i(p)$ is the unique element
of $A_i$ which is compatible with (and in fact weaker than) $ p$.
$h^j_i$ extends to a unique homomorphism
from $F_j$ to $F_i$, which we will also call $h^j_i$. 
\ENDent
\proc *commuting
 Fact: (1) The functions $h^j_i$ commute, i.e., if $i<j<k$ then 
$h^k_i= h^j_i\circ h^k_j$. 
\\
(2) If $i<j$, $p\in A_j$, then $p$ is stronger than $h^j_i(p)$.
\qed\commuting

Now let $\M =  (M, (F_i)_{i\in I}, (R_{i,z})_{i\in I, z\in F_i}, 
		(h^j_i)_{i\in I, j\in I, i<j})$, where we treat all
sets $F_i$, $R_{i,z}$, $h^j_i$ as relations on $M$.

\proc Definition: Let $G \subseteq P$ be a filter which meets
all the sets $A_i$, say $G \cut A_i = \{y_i(G)\}$.  Define $f_G:M\to M$
as follows: If $x\in F_i$, then $f_G(x) = x * y_i(G)$ (here $* = *_i$ is
the group operation on $F_i$).

\proc *fhasnofixpoints Fact:  If $G$ is a filter which meets all
sets $A_i$,  then $f_G$
is an automorphism of $\M$ without fixpoints. 

\proof: It is clear that the sets $F_i$ and the relations $R_{i,z}$
are preserved.  Note that for $i<j$ we have $h^j_i(y_j)=y_i$, since
$y_i$ and $y_j$ are compatible.  Since the functions $h^j_i$ are
homomorphisms, we have $h^j_i(f_G(x)) =
h^j_i(x*y_j)=h^j_i(x)*y_i=f_G(h^j_i(x))$, so also $h^j_i$ is
preserved.\QED\fhasnofixpoints

\smallskip

 So $\M$ is potentially nonrigid. 
So by  PA$(P,{\lambda})$ we know that $\M$ is really nonrigid.

Finally we will show how a nontrivial automorphism of $\M$ defines a
filter $G^*$ meeting all the sets  $A_i$. 

So let $F$ be an automorphism.  Let $1_i$ be the neutral element of
$F_i$, and assume $F(1_i) = w_i$.  Since the sets $F_i$ are predicates
in our structure, we must have $w_i\in F_i$.
Using the predicates $h^j_i$ we can  show:  If $i<j$, then
$h^j_i(w_j) = w_i$, and using the predicates
$R_{i,z}$ we can show that for all $z\in F_i$ we must have $F(z)
 = z * w_i$. 

Therefore, as $F$ is not the identity, 
 we can find $i^*\in I$ such that $w_{i^*}\not=
1_{i^*}$.  From now on we will work only with $I^ * = \{i\in I: i^* \le
i\}$. Since every antichain $A_i$ is refined by some antichain $A_j$
with 
$j\in I^*$ it is enough to find a directed set which meets all antichains
$A_j$ for $j\in I^*$. 

Let $u_i=\supp(w_i)$.  So for all $i\in I^*$ the set $u_i$ is finite
and nonempty (since $h^i_{i^*}[u_i] \supseteq u_{i^*} \not=\emptyset$).

\proc *lastfact Fact: 
\BEGINdent
\ite 
     1 If $J \subseteq I^*$ is a finite set, then there is a
family $\{p_i:i\in  J \}$ such that 
\itemitem{(a)} for all $i\in J$: $p_i\in u_i$
\itemitem{(b)} for all $i,j\in J$: If $i<j$, then $h^j_i(p_j)=p_i$.
\ite 
     2 There is a 
family $\{p_i:i\in  I^* \}$ such that 
\itemitem{(a)} for all $i\in I^*$: $p_i\in u_i$
\itemitem{(b)} for all $i,j\in I^*$: If $i<j$, then $h^j_i(p_j)=p_i$.
\ite 3 If $\{p_i:i\in  I^* \}$ is as in (2), then this set generates a
filter which will meet all sets $A_i$. 
\ENDent

\proof\ of (1): As $I^*$ is directed, we can find an upper
bound $j$ for $J$. Let $p$ be an  element of $w_j$ such that 
$p_i:= h^j_i(p_j)\in w_i$ for all $i\in J$. 

\noindent (2) follows from (1), by the compactness  theorem of
propositional calculus.  (Recall that all sets $u_i$ are finite.) 

\noindent (3): We have to show that for any 
$i_1, i_2 \in I^*$ the conditions  $p_{i_1}$ and $p_{i_2} $ are
compatible, i.e., have a common extension. Let 
$j$ be an upper bound of $i_1$ and $i_2$.  
Then $p_j$ witnesses that  $p_{i_1}$ and $p_{i_2}$ are compatible, as
$h^j_{i_1}(p_j)=p_{i_1} $ and 
$h^j_{i_2}(p_j)=p_{i_2} $.\hfill  \QED\lastfact\QED\PIisFFA. 

For the theorem \twotrees\ below we need the following definitions.
\proc Definition:   A tree on a set $X$ is a nonempty set $T$ of
finite sequences of elements of 
$X$ which is closed under restrictions, i.e., if
$\eta:k\to X$ is in  $T$ and $i<k$, then also $\eta \on i\in T$.   The
tree ordering $\le_T$ is given be the subset (or extension) relation: 
$ \eta \le \nu $ iff $ \eta \subseteq \nu $ iff $\exists i: \eta=\nu
\on i$. 

For $\eta\in T$ let $\Suc_T(\eta):= \{x\in X: \eta\concat x\in T \}$. 

For $A \subseteq T$, 
 $\eta\in T$ we let $\rk(\eta, A)$ be the rank of $\eta$
with 
respect to $A$, i.e., the rank of the (inverse) tree ordering 
on the set 
$$ \{ \nu: \eta \le \nu \in T, \forall \nu': \eta
\le \nu' <\nu \Limpl \nu' \notin A \}$$
In other words,  $\rk(\eta,A)=0 $ iff $\eta \in A$,
$\rk(\eta,A)=\infty$ iff there is an infinite branch of $T$ starting
at $\eta$ which avoids $A$, and $\rk(\eta,A) = \sup\{\rk(\nu,A)+1:  %\}
\nu$ a direct successor of $\eta %\{
\}$ otherwise.

\proc *twotrees  Theorem:  For any two structures $\A$ and $\B$ there
is a structure $\C=\C(\A,\B)$ such that in any extension $V' \supseteq
V$ of the universe, $V' \models `` \A \iso \B \  \liff \ \C $ is not
rigid.'' 

\proof: Wlog $|A| \le |B|$.   Also wlog $\A$ and $\B$ are structures
in a purely relational language $\L$, and we may also assume that
$A\cut B = \emptyset$.

We will say that a tree $ T $ on $A\cup B$ ``codes $A$'' iff
\BEGINdent
\ite 1 $\Suc_{T}(\eta) \in \{A, B \}$  for all $\eta\in T$. 
\ite 2 Letting $T^A:= \{\eta\in T: \Suc_{T}(\eta)= A\}$, the ranks
	$\rk(\eta, T^A)$ are $< \infty$ for all $\eta\in T$. 
\ite 3 The function $\eta\mapsto \rk(\eta, T^A \setminus \{\eta\})$ is 1-1 on $T^A$.  
\ENDent
Such a tree can be constructed inductively as $T = \bigcup_n T_n$,
where the $T_n$ are well-founded trees,  each $T_{n+1} $
end-extends $T_n$, and all  nodes in $T_{n+1} - T_n$ are from
$B$ except those at the top (i.e., those whose immediate successors
will be in $T_{n+2}-T_{n+1}$).   Because we have complete freedom in
what the rank of the  tree ordering for each connected component
of $T_{n+1} - T_n$ should be (and because all the $T_n$ have size
$=|B|$), we can arrange to satisfy (1), (2) and (3).  

Moreover, we can find trees $T_0$ and $T_1$, both coding $A$, such that 
\BEGINdent
\ite 4 
$\Suc_{T_0}(\emptyset)=A$,   $\Suc_{T_1}(\emptyset)=B$. 
\ENDent

We will replace  the roots ($\emptyset$) of the trees $T_0$ and $T_1$
by some new and distinct objects $\emptyset_0$ and $\emptyset_1$.  So
the trees $T_0 $ and $T_1$ will be disjoint (by (4)). 

Now define the structure $\C$ as follows:   We let $C = T_0 \cup T_1$.
%   Fortunately, $T_0$ and $T_1$ are
% already disjoint except for their roots, so  we will let 
% $C:= (T_0 \setminus\{\emptyset})\cup (T_1\setminus\{\emptyset\}) \cup
% \{\emptyset_0, \emptyset_1\}$, where $\emptyset_0$ and \emptyset_1$ are
% distinct (and not in $T_0\cup T_1$). 
% 
% T_1 \times \{1\}$).  

The
 underlying language of $\C$ will be the language $\L$ plus an
additional binary relation symbol $\le$, which is to be interpreted as
the tree order. 
Whenever $R$ is an $n$-ary relation in the language $\L$,
we interpret $R$ in $\C$ by 
 $$R^\C:= \{(\eta\concat a_1, \ldots, \eta\concat a_n):  %\}
  \vtop{\hsize=0.5\hsize\parindent0cm 
$\eta\in T_0 \cup T_1$, and \\
$\Suc(\eta)= A $ $\Limpl$ $(a_1 , \ldots, a_n)\in R^\A$, \\
$\Suc(\eta)= B $ $\Limpl$ $(a_1 , \ldots, a_n)\in R^\B$ $%\{
\}$}$$

Now work in any extension $V ' \supseteq V$.  First assume that $f:\A
\to \B$ is an isomorphism.  We will define a map
$g:T_0 \to T_1$ such that the map $g \cup g^{-1}$
is a (nontrivial)  automorphism of $\C$. 

$g$ is defined inductively as follows: 
\BEGINdent
\ite a 
      $g(\emptyset_0) = \emptyset_1$. 
\ite b If $\Suc_{T_0}(\eta)=\Suc_{T_1}(g(\eta))$, then $g(\eta\concat a)
= g(\eta) \concat a$. 
\ite c Otherwise, $g(\eta\concat a) = g(\eta) \concat f(a)$ or
 $g(\eta\concat a) = g(\eta) \concat f^{-1}(a)$, as appropriate. 
\ENDent

It is easy to see that $g\cup g^{-1}$ will then be a nontrivial automorphism. 

\smallskip
Now assume conversely that $g:\C\to \C$ is a nontrivial automorphism.
Recall that the tree ordering is a relation on the structure $\C$, so
it must be respected by $g$. 

First assume that there are
$i,j\in\{0,1\}$ and an $\eta$ 
such that
$$ (*)\ \ \ \  \eta\in T_i\ \ \ g(\eta)\in T_j \ \ \ 
\ \ \ \ \Suc_{T_i}(\eta)\not=\Suc_{T_j}(g(\eta))$$
   So without loss of
generality $\Suc_{T_i}(\eta)=A$ and $\Suc_{T_j}(g(\eta))=B$.  Now define a
map $f:A\to B$ by requiring 
$$ g(\eta \concat a ) =g(\eta) \concat f(a)$$
and check that $f$ must be an isomorphism. 

Now we show that we can always find $i,j,\eta$ as in  $(*)$. 
  If not, then we can first see that $g$ respects
$T_0$ and $T_1$, i.e., $g(\eta)\in T_0$ iff $\eta\in T_0$. 
  Next, our assumption implies that the
functions $g\on T_0$ respect the sets $T_0^A$, i.e., $\eta\in T_0^A $ iff
$g(\eta)\in T_0^A$.  Hence for all $\eta\in T_0^A$, $\rk(\eta, T_0^A) =
\rk(g(\eta), T_0^A)$, so (by condition (3) above)
 $g(\eta)=\eta$ for all $\eta\in T_0^A$.   Since
every $\nu\in T_0$ can be extended to some $\eta\in T_0^A$ and $g$
respects $<$, we must have $g(\nu)=\nu $ for all $\nu \in T_0$.   The
same argument shows that also $g\on T_1$ is the identity. \QED\thistheorem

\bigskip

\proc Remarks on other applications: Which other consequences of PFA
(see, e.g.,   [\PFApplications]) 
are already implied by BPFA?    On the one hand it is clear that if
PFA is only needed to produce a sufficiently generic function from
$\omega_1$ to $\omega_1$, then the same proof should show that BPFA is
a sufficient assumption.  For example: 
\BEGINdent
\item {} $\BPFA $ implies ``all $\aleph_1$-dense sets of reals are
isomorphic.''
\ENDent
%%%
% It is also possible to show that $\BPFA$ implies that the only gaps in
% $\P(\omega)/{\rm fin}$ where both sides have uncountable cofinality
% are the $(\omega_1, \omega_1^*)$ gaps, and that $\BPFA$ implies
% ${\frak b}= \aleph_2$. 
%%%

On the other hand, as we will see in the next section, the consistency
strength of $\BPFA$ is quite weak.  So $\BPFA $ cannot imply any
statement which needs large cardinals, such as ``there is an Aronszajn
tree on $\aleph_2$.''   In particular, BPFA does not imply PFA. 
      
We do not know if BPFA already decides the size of the continuum, but
Woodin has remarked that the bounded {\bf semi}proper forcing axiom
implies $2^{{\aleph_0}}=\aleph_2$. 

\neusection
\BEGINsection 2. The consistency of BPFA

\proc Definition:  For any cardinal $\chi$, $H(\chi)$ is the
collection of sets which are hereditarily of cardinality $< \chi$: 
Letting $trcl(x)$ be the
transitive closure of $x$, $trcl(x) = \{x\}\cup\bigcup x
\cup\bigcup\bigcup x \cup\cdots $, we have 
$$ H(\chi) = \{ x: |trcl(x) | < \chi \}$$
(Usually we require $\chi$ to be regular)

\proc Definition:  
Let $\kappa$ be an regular cardinal.   We say that $\kappa$ is
``reflecting'' or more precisely, ${\Sigma}_1$-reflecting, if: 
\\
For any first order formula $\varphi$ in  the language of set
theory, for any $a\in H(\kappa)$: 
\BEGINdent
\item{} IF there exists  a regular cardinal $\chi \ge  \kappa$ such that
$H(\chi)\models \varphi(a)$
\item {} THEN there is a cardinal ${\delta} <
\kappa$ such  
that $a\in H(\delta) $ and $H(\delta) \models \varphi(a)$. 
\ENDent

\proc *regrem Remark: (1)  We may require ${\delta}$ to be regular
without changing 
the concept of ``${\Sigma}_1$-reflecting''.
\\
(2) We can replace ``for all $\chi$'' by ``for unboundedly many $\chi$''

\proof: (1) Assume that $H(\chi) \models \varphi(a)$, $\chi$ regular.
Choose some large enough $\chi_1$ such that $H(\chi)\in H(\chi_1)$,
$\chi_1$ a successor cardinal.  So $H(\chi_1) \models ``\exists
\chi\,$, $\chi $ regular, $H(\chi)$ exists and $H(\chi) 
 \models `\varphi(a)$'\,''. We can
find a (successor) ${\delta}_1 < \kappa $ such that $H({\delta}_1)
\models ``\exists \delta\,$, $\delta $ regular, $H(\delta) \models
`\varphi(a)$'\,'' So ${\delta}$ is really regular. 

\noindent (2): If $\chi < \chi_1$ then 
$H(\chi) \models ``\varphi $'' iff $H(\chi_1) \models `` H(\chi)
 \models `\varphi $'\,''.  \QED\thistheorem

\proc *inacc Remark:  It is easy to see that if $\kappa$ is
reflecting,  then $\kappa$ is a strong limit, hence inaccessible.
Applying ${\Sigma}_1$ reflection, we get that $\kappa$ is
hyperinaccessible, etc.  \qed\inacc

\proc Remark: 
(1)   There is a closed unbounded class $C$ of cardinals such that
every regular $\kappa \in C$ (if there are any) is ${\Sigma}_1$
reflecting.  So if ``$\infty $ is Mahlo'', then there are many
$\Sigma_1$-reflecting cardinals. 
% 
% 
%   If ${\lambda}$ is a Mahlo cardinal, then $H(\lambda)
% \models ``$there is a reflecting  cardinal'' (in fact, there is a closed
% unbounded set $C \subseteq \kappa$ such that all regular cardinals in
% $C$ are reflecting in $H({\lambda})$.
 \\
 (2)
 If $\kappa$ is reflecting, then $L \models $ ``$\kappa$ is reflecting''.  

\proof: 
(1) For any set $a$ and any formula $\varphi$ let $f'(a,
\varphi) = \min\{\chi\in RCard: H(\chi) \models \varphi(a)\}$ (where
$RCard$ is the class of regular cardinals, and we define $\min\emptyset=0$).
Now let $f:RCard  \to RCard $ be defined by $f(\alpha) = \sup\{ f'(a,
\varphi): \hbox{ $\varphi$ a formula, $a\in H(\alpha)$}\}$, and let $C=
\{ {\delta}\in Card: \forall \alpha \in RCard \cut {\delta} \,\,
f(\alpha) < {\delta} \}$. 
% $C$ is a closed unbounded class of cardinals.  If $C$ contains a
% regular cardinal $\kappa$ then it is easy to see that $\kappa$ must be
% reflecting. 
% 
% Now if ${\lambda}$ is a Mahlo cardinal, then relativizing this
% argument to $H(\lambda)$ we get that every element in $Reg \cut
% C^{H(\lambda)}$ will be reflecting  in $H(\lambda)$. 
%
% 
% \proof\ of  (2): Clearly $\kappa$ will be strongly
% inaccessible in $L$.  Now 
% assume that $\chi$ is a cardinal in $L$,  $a\in L_\kappa$, and
% $L_\chi \models \varphi(a)$.  $\chi$ may not be a cardinal in $V$, but
% we can find a  true cardinal $\chi_1 > \chi$ in $V$.  
% So 
% $$ V_{\chi_1} \models ``\exists \chi\, L_\chi \models `\varphi(a) \txt'
% \hbox { and $\chi$ is a cardinal in $L$''} $$
% Since $\kappa $ was \reflecting L, we can find a cardinal
% ${\delta}_1 < \kappa $ such that 
% $$ V_{\delta_1} \models `` \exists {\delta} \, H(\delta) \models
% `\varphi(a)\txt'
% \hbox { and ${\delta}$ is a cardinal in $L$''} $$
% By the absoluteness of $L$ we have $H(\delta) ^{V_{\delta_1}} =
% H(\delta)$, so we get that $H(\delta) \models \varphi(a)$. 
% 
% Also by the absoluteness of $L$ we have $L^{V_{\delta_1}} =
% H(\delta_1) $. 
% So ${\delta} $ is a cardinal in $L_{{\delta}_1}$, but as ${\delta}_1$
% is itself a cardinal in $L$ we get that also ${\delta}$ is a cardinal
% in $L$.

\noindent
(2) is also  easy. 
  \qed\thistheorem

Our main interest in this concept is its relativization to $L$.  In
this context we recall the following fact: 

\proc *LisH Fact:  Assume $V=L$. Then for all (regular) cardinals
$\chi$, $H(\chi) = L_\chi$. \qed\thistheorem

\proc *Hchi Fact:  Assume $P\in H({\lambda})$ is a forcing notion,
$\chi > 2^{2 ^{\lambda} }$ is regular. Then 
\BEGINdent
\ite 1 For any $P$-name $\name x$ there is a $P$-name $\name y\in
H(\chi)$ such that $\forces_P ``\name x\in H(\chi) \Limpl \name x =
\name y $''.   (And conversely, if $\name x\in H(\chi)$, then
$\forces_P ``\name x \in H(\chi)$''.)
\ite 2 If $\name x \in H(\chi)$, $\varphi(\cdot)$ a formula, then $$
\forces ``H(\chi) \models \varphi(\name x) \txt'' \ \Liff \ 
 ``H(\chi) \models\, ` \forces\varphi(\name x)\txt'\, \txt''$$
\ENDent
Proof: (1) is by induction on the rank of $\name x$ in $V^P$,
 and (2) uses (1). \qed\thistheorem

\proc *properisabsolute Fact:  Let $P$ be a forcing notion, $P\in
H({\lambda})$, $\chi > 2^{2 ^{\lambda}}$ regular.  Then $P$ is proper
iff $H(\chi) \models ``P$ is proper''.  \qed\thistheorem

\proc *findsmallproper Lemma:  
Assume that  $\kappa$ is reflecting,  ${\lambda} < \kappa$ is a regular 
cardinal, $\A$ and $\B$ are structures in $H({\lambda})$. 
\\
If there is a proper forcing notion $P$ such that $\forces_P ``
\A \iso \B$'', then there is such a 
(proper) forcing notion in $H(\kappa)$. 

\proof:  Fix $P$, and let $\chi$  be a large enough regular cardinal.
So  $H(\chi) \models ``P $ proper, $P\in H(\mu)$,
$\bigl(2^{2^\mu}\bigr)$ exists''.  
Also, there is a $P$-name $\name f\in H(\chi)$ such that 
  $\forces_P ``H(\chi) \models `\name f: \A \to \B$ is an isomorphism'\,'', 
 so by \Hchi(2),  
  $H(\chi) \models `` \forces_P `\name f: \A \to \B$ is an
isomorphism'\,''. 

Now 
we use the fact that $\kappa$ is reflecting.  We can find
${\delta} < \kappa$, ${\delta} > {\lambda}$
 and $\chi' \in H({\delta})$ such that
$H({\delta}) \models `` \exists \nu\,  \exists Q\in H(\nu), Q$ proper,
$\exists  \name g \,  
 \forces_{Q} `\name g:\A \to \B$ is an isomorphism', and 
$\bigl(2^{2 ^{\nu}}\bigr)$ exists.''
So this $Q$ is really proper, and $Q$ forces that $\A$ and $\B$ are
isomorphic. 
 \QED\thistheorem

\proc *smallforcing Fact:  If $\kappa$ is reflecting, 
$P\in H(\kappa)$
is a forcing notion, then $\forces_P ``\kappa$ is reflecting''. 

\proof:   Let $P\in H({\lambda})$, ${\lambda} < \kappa$. 
Assume that $p \forces `` H( \chi) \models `\varphi(\name
a)\txt', \name a\in H(\kappa)$''.  We may assume that $\name
a\in H(\kappa)$.  By \Hchi\  we have 
$H(\chi) \models ``p \forces
`\varphi(\name a )$'\,'', so there is a ${\delta} < \kappa$,
${\delta}> {\lambda} $, 
such that $H({\delta}) \models  ``p \forces
`\varphi(\name a )$'\,'', hence 
 $p \forces `` H({\delta} ) \models `\varphi(\name a)$'\,''.
${\delta}$ is a cardinal in $V^P$, because $|P| < {\lambda} <
{\delta}$. 
\QED\smallforcing

\proc Theorem: If ``there is a reflecting  cardinal'' is consistent
with ZFC, then also  
PI$({\rm proper}) $ (and hence BPFA, by \PIisFFA) is consistent with ZFC.

\proof: (Short version) 	We will use an CS iteration of length
$\kappa$, where $\kappa$ reflects.  All intermediate
forcing notions will have hereditary size  $< \kappa$. 
 By a bookkeeping argument we can take care of all possible 
pairs of structures on $\omega_1$. If in the intermediate model there
is a proper forcing notion making two structures isomorphic, then
there is such a forcing notion of size $< \kappa$, so we continue.
Note that once two structures have been made isomorphic, they continue
to stay isomorphic. 

\proof: (More detailed version) 
Assume that $\kappa$  reflects. We define a countable
support iteration $(P_i, Q_i:i < \kappa)$ of proper forcing
notions and a sequence $\< \name \M_i, \name \N_i: i <
\kappa> $ with 
the following properties for all $i < \kappa$: 
\BEGINdent
\ite 1 $P_i \in H(\kappa) $ 
\ite 2 $Q_i$ is a $P_i$-name,  $\forces_{P_i} $ ``$Q_i$ is proper,
	$Q_i \in H(\kappa) $''.  
\ite 3 $\forces_{P_i} 2^{\aleph_1} < \kappa $.   (This follows from
(1) and (2))
\ite 4 $\name \M_i$ and $\name \N_i$ are names for structures on
	$\omega_1$. 
\ite 5 $\forces_{P_i}$ ``If $\name\M_i \mathop\iso\limits_{\rm
	proper,\,{<}\kappa} \name\N_i$, then $\forces_{Q_i} `\name\M_i \iso
	\name\N_i$'\,''. 
\ENDent
With the usual bookkeeping argument we can also ensure that
\BEGINdent
\ite 6 Whenever  $\name \M$ and $\name \N$ are $P_i$-names for
	structures on $\omega_1$ for some $i$, then there are
	unboundedly (or even stationarily) many $j>i$ with 
   $\forces_j$ ``$\name \M_j=\name \M$,  $\name \N_j=\name \N$''
\ENDent
From (1) we also get the following two properties: 
\BEGINdent
\ite 7 $P_\kappa \models \kappa$-cc
\ite 8 Whenever $\name \M$ is a  $P_\kappa$-name for a
	structure on $\omega_1$, then there are $i< \kappa $ and
	a $P_i$-name $\name \M'$ such that $\forces_{\kappa} \name \M =
	\name \M'$. 
\ENDent
From these properties we can now show $\forces_\kappa \BPFA$.
$P_\kappa$ is proper, so $\omega_1$ is not collapsed. 
 Let
$p$ be a condition, and let $\name \M$ and $\name \N$ be
$P_\kappa$-names for structures on $\omega_1$, and assume that 
$$  p \forces_\kappa  \hbox{``$\name Q$ proper, } \forces_{\name Q} \name \M
	\iso \name \N \hbox{''}$$
where $\name Q$ is a $P_\kappa$-name. 
So by (8) we may assume that for some large enough $i < \kappa $
$\name \M $ and $\name \N$ are $P_i$-names. By (6) wlog we may assume
that $\name \M = \name \M_i$, $\name \N = \name \N_i$. 
Now letting $R$ be the
$P_i$-name $(P_\kappa / G_i) * \name Q$, we get 
$$ p \forces_i ``\forces_R \name \M \iso \name \N\txt''$$
But by \smallforcing,  $\forces_i$ ``$\kappa$ is reflecting'', so by 
the definition of $Q_i$ and by \findsmallproper\ we get 
that $p \forces _{i+1} \name \M\iso \name \N$. 
\QED\thistheorem

\proc Remark: Since \properisabsolute\ is also true with ``proper''
replaced by ``semiproper'', we similarly get that the consistency of a
${\Sigma}_1$-reflecting cardinal implies the consistency of the
bounded semiproper forcing axiom. \QED\thistheorem

\neusection
\BEGINsection 3. Sealing the $\omega_1$-branches of a tree

In this section we will define  a forcing notion 
%(originally from
%[???], see also chapter XI of [\ProperForcing] and of
%[\ProperImproper]) 
which makes the set
 of branches of an $\omega_1$-tree absolute.

\proc Definition: Let $T$ be a tree of height $\omega_1$.
  We say that $B \subseteq T$ is
an $\omega_1 $-branch if $B$ is a maximal linearly ordered subset of
$T$ and has  order type $\omega_1$.

\proc *whatisPT Lemma:  Let $T $ be a tree of height $\omega_1$.
  Assume that every node of $T$ is on
some $\omega_1$-branch, and that there are at uncountably many
$\omega_1$-branches.  (These assumptions are just to simplify the
notation).   Then there is a proper forcing notion $P'_T$ (in fact,
$P'_T$ 
is a composition of finitely many ${\sigma}$-closed and ccc forcing
notions)  forcing the following: 
\BEGINdent
\ite 1 $T$ has $\aleph_1$ many
	$\omega_1$-branches, i.e., there is a function $b:\omega_1
	\times \omega_1 \to T$ such that each set $B_\alpha
	=\{b(\alpha,{\beta}): {\beta}  <
	\omega_1\}$ is a an end segment of of a branch of  $T$
	(enumerated in its natural 
	order), and every $\omega_1$-branch is (modulo a countable
	set) equal to one of the $B_\alpha $s, and the sets $B_\alpha$ are
	pairwise	disjoint.   
\ite 2  There is a function $g:T\to \omega$ such that for all $s < t $
	in $T$, if $g(s)=g(t) $ then there is some (unique) $\alpha  <
	\omega_1$ such 
	that $\{s,t\} \subseteq B_\alpha $. 
\ENDent

The proof consists of two parts.  In the first part (\mitchelfact) we
show that we may wlog assume that $T$ has actly $\aleph_1$ many
branches.  This observation
 is a special case of a theorem of Mitchell
[\Mitchell, 3.1].\\
 In the second part we describe the forcing notion which works under
the additional assumption that $T$ has only $\aleph_1$ many branches.
This forcing notion is esentially the same as the one used by
Baumgartner
in  [\IteratedForcing, section 8].

\proc *mitchelfact Fact: Let $T$ be a tree of height $\omega_1$,
$\kappa>|T|$, and let $R_1$ be the forcing notion adding $\kappa$ many
Cohen reals.  In $V^{R_1}$, let $R_2$ be a ${\sigma}$-closed focing
notion.   Then every branch of $T$ in $V^{R_1*R_2}$ is already in
$V^{R_1}$ (and in fact already in $V$). \\
Hence, taking $R_2$ to be the Levy collapse of the number of branches
of $T$  to $\aleph_1 $ (with countable conditions), $T$ will have at most
$\aleph_1$ many branches in $V^{R_1*R_2}$. 

\proof : Assume that $\name b$ is a name of a new branch. So the set 
$$T_{\name b}:=\{t\in T: \exists p\in R_2\, p \forces t\in \name b\}$$
is (in $V^{R_1}$) a perfect subtree of $T$. In particular, there is an
order-preserving function $f:2^{<\omega} \to T_{\name b}$.  Since
$\kappa$ was chosen big enough, we can find a real $c\in 2^\omega \cap
V^{R_1}$ which is not in $V[f]$.  Now note that $T'$ is
${\sigma}$-closed, so there is $t^*\in T$ such that 
$ \forall n\, f(c\on n) \le t^*$.  But this implies that 
$$c = \bigcup \{ s\in 2^{<\omega}: f(s) \le t^*\}$$
can be computed $V[f]$, a contradiction. \QED\thistheorem

Now we  describe a forcing notion $P'_T$ which works under the
assumption that  $T$ has
not more than $\aleph_1$ branches.  In the general case we can then
use the forcing $P_T = R_1 * R_2 * P'_T$. 

\proc  Definition:
Let $T$ be a tree of height $\omega_1$ with $\aleph_1$ many
$\omega_1$-branches
$\{B_i:i<\omega_1\}$ and assume that each note of $T$ is on some
$\omega_1$-branch.
Let $B'_j=B_j\setminus \bigcup_{i<j}B_i$,
$x_j=\min (B_j')$ so that the sets $B'_j$ are disjoint end segments of
the branches $B_j$, and they form a partition of $T$.
Let $A=\{x_i:i<\omega_1\}$.
% (so $A$ does not include any linearly ordered uncountable set)

The forcing ``sealing the branches of $T$'' is defined as
$$P'_T=\{f:f \hbox{ a finite function from $A$ to $\omega$,
and if $x<y$ are in $\dom(f)$, then $f(x)\not=f(y)$}\}$$

\proc *cccofPT Lemma:  $P'_T$ satisfies the countable chain condition.
(In fact, much more is true: If $\langle p_i:i<\omega_1\rangle$ are
 conditions  in $P$, then there are  uncountable sets $S_1$, 
$S_2\subseteq \omega_1$ such that whenever
$i\in S_1$,
$j\in S_2$, then $p_i$ and $p_j$ are compatible.  See
[\ProperImproper, XI]) 

\proof: Essentially the same as in [\IteratedForcing, 8.2].
\qed\thistheorem

To conclude the proof of \whatisPT, note that any generic filter $G$
on $P_T'$ induces a generic $f_G:A \to \omega$. 
  Let $g_G:T \to \omega $ be defined by $g_G(y)=f_G(x_i)$ for all $y\in
 B_i$.   This function $g_G$ fulfills the requirement \whatisPT(2).
 \QED\whatisPT

\neusection
\BEGINsection 4. BPFA and reflecting cardinals  are equiconsistent

In this section we will show that 

\proc *equicon Theorem: 
If \BPFA\ holds, then the cardinal $\aleph_2$ (computed in $V$) is  
${\Sigma}_1$-reflecting in $L$. 

Before we start the proof of this theorem, we show some 
general properties of ``sufficiently generic'' filters. 

First a remark on terminology: When we consider $\BFA(P,\lambda)$, then
by ``for all sufficiently generic $G^*
\subseteq P$, $\varphi(G^*)$ 
holds'' we mean:  ``there is a $P$-name $\name f: \lambda \to \lambda $
 such that:  whenever a filter $G^*$ interprets $\name f$, then
$\varphi(G^*)$ will hold''.   A description of the name $\name f$ can
always be deduced from the context.  Instead of a single name  $\name
f$ we usually have a family 
of $\lambda  $ many names.

The first lemma shows that from any sufficiently 
generic filter we can correctly
compute the first order theory (that is, the part of
it which is forced), or equivalently, the first order
diagram, of any small structure in the extension.

\proc *absolute Lemma: Let $P$ be a forcing notion, $\forces_P ``\M$
is a structure with universe $\lambda$ with $\lambda$ many relations
$(\name R_i: i< \lambda)$''.  Assume $\BFA(P,\lambda)$.  Then for every
sufficiently generic filter $G^* \subseteq P$, 
letting $\M^* =(\lambda, \name R_i[G^*])_{i< \lambda}$,  (where $\name
R_i[ G^* ]:= \{(\VEC x)\in {\lambda}^n: \exists p\in G^*\, p \forces 
\M \models \name R_i(\VEC x)\}$)   we have: 
\BEGINdent
\item{}  Whenever $\varphi$ is a closed formula such that $\forces_P
	\M \models \varphi$, 
\item{} then $\M^* \models \varphi$. 
\ENDent
\proof:    Let $\chi$ be a large enough cardinal, and
 let $N$ be an elementary submodel of $H(\chi)$ of size $\lambda$
 containing all the necessary information (i.e., $\lambda \subseteq N$,
 $(P, \le)\in N$, $(\name R_i:i< \lambda)\in N$).  

By $\BFA(P,\lambda)$ we can find a filter $G^* \subseteq P$ which
 decides all $P$-names of elements of $\M$ which are in $N$ and all
 first order statements about $\M$,  i.e., 
\BEGINdent
\ite 1 For 
 all $\name \alpha \in N$, if $\forces_P ``\name \alpha \in \lambda$''
 then there
 is ${\beta}\in \lambda$ and $p\in G^*$ such that $p \forces_P `` \name
 \alpha = \check {\beta}$''.  
\ite 2 For all $\VEC{ \alpha }\in {\lambda} $ and all formulas
 $\varphi(\VEC x)$  there is $p\in G^*$ such that either $p
 \forces ``\M \models \varphi(\VEC{ \alpha })$'' or  $p \forces ``\M \models
 \lnot \varphi(\VEC{\alpha })$''. 
\ENDent

We now claim that for every formula $\varphi(\VEC x)$, for every $\VEC
 {\name a}\in N$:  If $\forces_P ``\M \models \varphi(\VEC {\name a})$'',
 then $\M^* \models \varphi(\name a_1[G^*], \ldots, \name a_k[G^*])$.
We assume that $\varphi$ is in prefix form, so in particular negation
signs appear only before atomic formulas. 
The proof is by  induction on the complexity of $\varphi$, starting
from atomic and negated atomic formulas.   We will only treat the case
$\varphi = \exists x \, \varphi_1$.   So assume 
% For atomic and negated atomic formulas this follows immediately
% follows from (1) and (2).The cases $\varphi= \varphi_1 \logand
% \varphi_2$ and $\varphi =\forall x \, \varphi_1(x)$ are trivial.  In
% the case $\varphi=\varphi_1 \lor \varphi_2$ we use    
that  $\forces_P \M \models \exists x \varphi_1(x, \VEC {\name a})$. 
We can find a name $\name b\in N$ such that $\forces_P \M \models 
\varphi_1(\name b, \VEC {\name a})$, so by induction hypothesis we get
$ \M^* \models \varphi_1(\name b [G^*], \name a_1[G^*], \ldots, \name
a_n[G^*])$.  \QED\absolute

\proc Remark:  In a sense the previous lemma characterizes
``sufficiently generic'' filters.  More precisely, the following is
(trivially) true:  Let $P$ be a complete Boolean algebra, let
$\forces_P \name f: \lambda \to \lambda$, and let $\name \M = (\lambda,
\name f)$, where we treat $\name f$ as a relation.  For any
ultrafilter $G^* \subseteq P$ the model $\M^* = (\lambda, \name
f[G^*])$ is well-defined.  Since $\name f$ is forced to be a function,
we have $\forces_P ``\M \models `\forall 
\alpha \, \exists {\beta}\,\, (\alpha, {\beta})\in \name f$'\,''.
Clearly $G^*$ ``decides''  $\name f$ (as a function)  iff 
 $\M^*$ satisfies the same $\forall \exists$ statement. \QED\thistheorem

This last remark  suggests the following easy characterization of
$\BFA(P)$: 
\proc *howtodefineit Definition:  Let $P$ be an arbitrary forcing
notion, not necessarily a complete Boolean algebra.
 If $\name f$ is a $P$-name of a function from ${\lambda}$ to
${\lambda}$, then let the ``(forced) diagram'' 
 of $\name\M = (\lambda, \name f) $ be defined by 
$$D^{\scriptscriptstyle\forces}(\name M)=
 D^{\scriptscriptstyle\forces}(\name f)= 
\{(\varphi, \alpha_1, \ldots, \alpha_n):    %\}
  \vtop{\hsize=0.5\hsize\parindent0cm
	$\varphi(x_1, \ldots , x_n)$ a first order formula,
$\alpha_1$, \dots, $\alpha_n\in {\lambda}$, $
\forces_P \varphi(\alpha_1, \ldots, \alpha_n)$ $%\{
	\}$}$$
The ``open (forced) diagram'' $D_{qf}^{\scriptscriptstyle\forces}(\name f)$ 
is defined similarly, but $\varphi$
ranges only over quantifier-free formulas. 

\proc *freeofro Definition:  For any forcing notion $P$ let 
$\BFA'(P,\lambda)$ be the statement
$$ \BFA'(P,\lambda) = \vtop{\hsize0.7\hsize Whenever $\name f:
{\lambda} \to {\lambda}$ is a $P$-name of a function, then there is a
function $f^*$ such that $(\lambda , f^* ) \models
D_{qf}^{\scriptscriptstyle\forces}(\name f)$.}$$ 

\proc *secondCBA Fact:  For any forcing notion $P$, 
  $\BFA(P,{\lambda}) $ iff $\BFA'(P,{\lambda})$.

\proof:  $\BFA'(P,{\lambda})$ is clearly equivalent to
$\BFA'({\rm ro}(P),{\lambda})$.   The same is true (by definition) for $\BFA$.
So we may  wlog assume that $P$ is a complete Boolean Algebra. It is clear
that $\BFA(P,{\lambda}) \Limpl \BFA'(P,{\lambda})$. 
\\
Conversely, if 
$f^*$ is a function as in $\BFA'$, then the we claim that 
the  set $\{[\![ \name f(\alpha) = f ^* (\alpha) ]\!]: \alpha \in
{\lambda}\}$ generates a 
filter on $P$ (where $[\![\varphi ]\!]$ denotes the Boolean value of a
closed statement $\varphi$).   Proof of this claim:  If not, then
there are  ordinals $\alpha_1$, \dots, $\alpha_n$, ${\beta}_1$, \dots,
${\beta}_n$ such that $$f^*(\alpha_1)= {\beta}_1 \logand \cdots
\logand f^*(\alpha_n) = {\beta}_n$$ but the Boolean
value 
$$ [\![	
\name f(\alpha_1) = {\beta}_1 \logand 
\cdots \logand \name f(\alpha_n)={\beta}_n 
 ]\!]$$
is $0$.   This is a contradiction to the fact that $f^*$ witnesses 
$\BFA'(P,\lambda)$. 
\QED\thistheorem

After this digression we now continue our preparatory work for the
proof of theorem \equicon.  Our next lemma shows that a generic filter
will  not only reflect first order statements about small structures,
but will also preserve their wellfoundedness. 

\proc *whyWF Lemma:
Assume that  $\forces_P ``\name \M  = (\lambda,\name  E)$ is a
 well-founded structure, $\lambda$ is a cardinal''.  Assume that
 $\cf(\lambda) > \omega$, and assume that $\BFA(P,\lambda)$ holds.  Then
 for every sufficiently generic filter $G^* \subseteq P$  we
 have that $\M^* := (\lambda, \name E[G^*])$ is well-founded. 

(We will use this lemma only for the case where $P$ is proper and
 $\lambda = \omega_1$.)
\proof:  For each $\alpha < \lambda$ let  $\name r_\alpha$ be the name 
 of the canonical  rank function on $(\alpha, \name E)$, i.e., 
$$ \forces_P ``\dom(\name r_\alpha) = \alpha, \ 
\forall {\beta} < \alpha \,\, \name r_\alpha ({\beta})= \sup
	\{\name r_\alpha({\gamma})+1: {\gamma} E {\beta} \}\hbox{''}$$
As  $\forces_P `` {\lambda} $ is a cardinal'', 
 we have  $\forces_P ``\range(\name r_\alpha) \subseteq
 \lambda$'', so any sufficiently generic filter $G^*$ will interpret
 all the functions $\name r_\alpha$.   Applying lemma \absolute\ to
 the structure $(\alpha , \name E[G^*], \name r_\alpha[G^*])$ we see
 that $\name r_\alpha[G^*]$ is indeed a rank function witnessing that
 $(\alpha, \name E[G^*])$ is  well-founded.  Since $\cf(\lambda) >
 \omega$ this now implies that also $(\lambda, \name E[G^*])$ is
 well-founded. \QED\whyWF

We now start the proof of \equicon.  The definitions in the following
paragraphs will be valid throughout this section. 

Assume $\BPFA$.  Let $\kappa:=\aleph_2$. We will show that
$\kappa$ is reflecting in $L$. 
  It is clear that $\kappa$ is regular in $L$.

\proc *coveringlemma Claim:  Without loss of generality we may assume:
\BEGINdent
\ite 1 $0^\#$ does not exist, i.e., the covering lemma holds for $L$.
\ite 2 $\aleph_2^{\aleph_1} =  \aleph_2$.
\ite 3 There is $A \subseteq \aleph_2 $ such that 
        whenever $X \subseteq Ord$ is
	of size $\le \aleph_1 $, then $X\in L[A]$. 
\ENDent
\proof: (1) If $0^\#$ exists, then $L_\kappa  \prec  L$,  and it is
	easy to see that this implies that $\kappa $ is a
	reflecting  cardinal in~$L$. 

(2) Let $P=\Levy(\aleph_2, \aleph_2^{\aleph_1} )$, i.e., members of $P$
	are partial functions from $\aleph_2$ to $\aleph_2^{\aleph_1}$
	with bounded domain. Since $P$ does not add new subsets of
	$\aleph_1$ and $P$ is proper, also $V^P$ will satisfy
	PI(proper,$\aleph_1$). 
	Also $\aleph_2^V = \aleph_2 ^{V^P}$ and $V^P \models
	\aleph_2^{\aleph_1 } = \aleph_2$, so we can wlog work in $V^P$
	instead of $V$. 

(3) By (2) we can find a set $A \subseteq \aleph_2$ such that
$\aleph_2^{L[A]} = \aleph_2$ and  every
	function from $\aleph_1$ to $\aleph_2$ is already in $L[A]$.
	By (1), every set $X$ of ordinals 
	of size $\le \aleph_1$ can be covered by a set $Y\in L$, $|Y|
	= \aleph_1$. Let $j:Y \to otp(Y)$ be order preseving, then
	$j[X] \in L[A]$, $j\in L$, so $X\in L[A]$. 
\QED\coveringlemma

\proof\  of \equicon: 
Let $\varphi(x)$ be a formula, $a\in L_\kappa$, and assume that $\chi
> \kappa$, $L_\chi \models \varphi(a)$, $\chi $ a regular cardinal in $L$.
We have to find an $L$-cardinal $\chi'< \kappa$ such that $a\in
L_{\chi'}$ and $L_{\chi'} \models \varphi(a)$. 

By \regrem, we may assume that $\chi$ is a cardinal in $L[A]$ or even
 in $V$. 

Informal outline  of the proof: We will define a forcing notion $P$.
In $V^P$ we will construct a model 
 $\M = (M,\in, \chi, x, \ldots ) \prec V^P $ of
size $\aleph_1$ containing all necessary information.  This model has an
isomorphic copy $\bar \M $ with underlying set $\omega_1$.  We will find a
``sufficiently generic'' filter $G^*$ which will ``interpret'' $\bar \M$
as $\M^*$.  By \whyWF\  we may assume that 
 $\M^* = (\omega_1, E^*, \chi^*\ldots) $ will be
well-founded, so we can form its transitive collapse 
$\M' = (M', \in, \chi',\ldots)$. By \absolute\ 
we have that $\M'
\models ``  \chi' $ is a cardinal in $L$'', i.e., $\chi'$ is a
 cardinal in $L_{M' \cut Ord }$.   The main point will be to show that
any filter on
 our forcing notion $P$ will code enough information to enable us to
 conclude that $\chi' $ is really a cardinal of $L$.

\proc Definition of the forcing notions $Q_0$ and $Q_1$: 
Let $Q_0$ be the Levy-collapse of $L_\chi[A]$ to $\aleph_1$, i.e.
the set of countable partial functions from $\omega_1$ to
$L_\chi[A] $ ordered by extension.

In $V^{Q_0}$ let $T$ be the following tree:  Elements of $T$ are
of 
the form $$( \< \mu_i: i < \alpha >, \< f_{ij}: i \le j < \alpha>)$$ (we
 will usually write them as 
$ \< \mu_i, f_{i j} : i \le j < \alpha>$), where the $\mu_i$
are ordinals $<\chi $, the $f_{ij}$ are a system of commuting
order-preserving embeddings, and $\alpha < \omega_1$.  $T$ is ordered
by the relation ``is an initial segment of''. 

If $B$ is a branch of $T$ (in $V^{Q_0}$, or in any bigger universe) of
length ${\delta}$ then $B$ defines a directed system $\<\mu_i, f_{ij}:
i\le j < {\delta}>$ of well-orders.  We will call the direct limit of
this system $(\gamma_B, <_B)$.  In general this may not be a well-order,
but it is clear that if the length of $B$ is $\omega_1$, then
$(\gamma_B, <_B)$ will be a well-order. 

Let $Q_1= P_T$ be the forcing ``sealing the $\omega_1$-branches of $T$''
described in \whatisPT.  We let $P = Q_0 * Q_1$. So $P$ is a proper
order, in fact it is a finite iteration of ${\sigma}$-closed and ccc
partial orderings.

\proc Definition:  In $V^P$ we define a model $\M$ as follows:
Let $\Omega$ a large enough regular cardinal of $V$. 
Let $(M,\in)$  be an elementary submodel of $(H(\Omega)^{V^P}, \in)$
of size $\aleph_1$
containing all necessary information, in particular $M \supseteq
L_ \chi[A]$.   We now expand $(M,\in)$ to a model $\M = (M,\in,
\chi, A, \ldots)$ by adding the following functions,
relations and constants: 
\BEGINdent
\item - a constant for each element of $L_{\xi} $ (where ${\xi}$ is
chosen such that $a\in L_{\xi}$)
\item - relations $M_0$ and $M_1$ which are interpreted as $M\cut
H(\Omega)^V$ and $M\cut H(\Omega)^{V^{Q_0}}$, respectively.
\item - constants $\chi$, $A$, $\kappa$,  $T$, $g$, $b$ ($b$ is the  function
enumerating the branches of $T$ from \whatisPT,  and $g$ is the specializing
function $g:T \to \omega$ also from \whatisPT)
\item - a function $c: \chi \times \omega_1 \to \chi$ such that for
	all ${\delta} < \chi$:  If $ \cf({\delta}) = \aleph_1$, then
	$c({\delta}, \cdot): \omega_1 \to {\delta} $ is increasing and
	cofinal in 	${\delta}$. 
\ENDent

Since $M$,  the underlying set of $\M$, is of
cardinality $\aleph_1$, we can find an isomorphic model 
$$\bar \M = (\omega_1, \bar E, \bar\chi, \ldots )$$

In $V$ we have names for all the above: $\name{\bar\M}$, 
 $\name{\bar E}$, etc. Now let $G^*$ be a sufficiently generic
  filter, i.e., $G^*$ will interpret all these names.  Writing
 $E^*$ for $\name {\bar E}[G^*]$, etc., and letting $\M^* = (\omega_1,
 E^*, \chi^*, \ldots)$, we may by \whyWF\ and \absolute\ assume that
 the following holds:

\proc *WFandTheory Fact: 
\BEGINdent
\ite 1 $(\omega_1, E^*)$ is well-founded. 
\ite 2 If $\psi$ is a closed formula such that $\forces_P ``
\M \models
\psi$'', then $\M^* \models \psi$. 
\ENDent

\proc Main definition:  We let 
$$\M' = (M', \in, \barchi, \ldots )$$
be  the Mostowski collapse of $\M^*$.  This is possible by
 \WFandTheory(1).    
%Let ${\eta} $ be the
% ``height''
% of $\M'$, i.e., ${\eta} = \M' \cut Ord $. 
$\M'_0 = (M'_0, \in)$ and $\M_1' = (M_1', \in)$ will be ``inner
models'' of $\M'$.

{\it 
\noindent {\bf Note:} We will now do several computations and
 absoluteness arguments 
 involving the universes $V$, $L[A']$, $\M'$,
 $L[A']^{\M'}=L_{M'\cut Ord}[A']$, etc.    By default, all set-theoretic
 functions, quantifiers, etc., are to be interpreted in $V$, but we
 will often also have to consider relativized notions, like
  $\M' \models ``L[A'] \models `\dots$'\,'' (which is
 of course equivalent to $L_{M'\cut Ord}[A] \models `\dots$'), or
 $\cf^{L[A']}$, etc. 
 }

%We let $\M_0$ and $\M_1$ be the restriction of the model $\M$ to $M_0$
% and $M_1$, respectively (of course some relations and functions, like
% $T$ or $g$ are not well-defined).  We define $\bar \M_0$, $\M_0^*$,
% $\M_0'$, $\bar \M_1$, etc. similarly. 
Note that $\M' \models `` L[A'] \models `\kappa ' = \aleph_2
\txt'\,\txt''$, so we get $\aleph_1^{\M'} = \aleph_1 ^V$. 

 We will finish  the proof
of \equicon\  with the following two lemmas: 

\proc *thinksphi Lemma: 
$a\in L_\barchi$, $L_{\barchi } \subseteq M'$ and $L_\barchi
\models \varphi(a)$.

\proc *isacardinal Lemma: $L \models \barchi$ is a cardinal. 
      
\proof\  of \thinksphi:  Since $\chi'+1 \subseteq M'$ and $\M'$
 satisfies a large fragment of ZFC, we have $L_\barchi  \subseteq M '
 $ and $L_\barchi \in M'$.  For each $y\in L_{\xi} $ let $c_y$ be the
 associated constant symbol, then by induction (using
 \WFandTheory(2)) it is easy to show that $y=c_y ^{\M'}$ for all $y\in
 L_{\xi}$.  Since $\forces_P ``\M \models
 \left[ L_\chi \models ` \varphi(a)\txt'\right]$'', we
 thus have $\M' \models ``L_\barchi  \models`
 \varphi(a)\txt'\,\txt''$.  But 
 $L_\barchi \subseteq M'$, so $L_\barchi \models \varphi(a)$. 
\QED\thinksphi

So we are left with proving \isacardinal. 
In $L[A']$ let $\mu$ be the cardinality of $\chi'$, and (again in
$L[A']$) let $\nu $ be the successor of $\mu$. 
We will prove \isacardinal\ by showing the following fact: 

\proc *nuisthere Lemma:  
 $\nu \subseteq M'$.

\proof\ of \isacardinal (using \nuisthere): In fact we show that
\nuisthere\ implies that $\barchi$ is a cardinal even in $L[A']$:
If not, then $\mu < \chi'$, and since $\nu $ is a cardinal in
 $L[A']$ we can find a ${\gamma}< \nu$ such that 
$L_{\gamma} [A' ] \models $``there is a function from $\mu $ onto
$\barchi $''.  By \nuisthere, $\gamma \in M'$, so by the well-known
absoluteness properties of $L$ we have $L_{\gamma}[A'] \subseteq M'$,
so $\M' \models `` L [A' ] \models ` \barchi $ is not a
cardinal.'$\,$'' But we also have 
$ \forces_P \M \models `` L [ A ] \models ` \chi $ IS a
 cardinal'\,'',
so we get a contradiction to \WFandTheory(2). \QED\isacardinal

\proof\ of \nuisthere:  We will distinguish two cases, according to
what the cofinality of $\mu$ is.   
% (Magidor has pointed out that case
% 1 can actually be treated in the same way as case 2 \dots 

\noindent{\bf Case 1:} $\cf(\mu)={\aleph_0}$. 
(This is the ``easy''  case, for which we do not need to know anything
 about the forcing $Q_1$ other than that it is proper, so the class
$\{{\delta}:\cf({\delta})={\aleph_0}\}$ is the same in $V$, $V^ {Q_0}$,
$V^P$, $L[A]$).
  We start our investigation of case 1 with
 the following remark: 

\proc Fact:  
\BEGINdent
\ite 1 For all ${\delta}$: If $\cf^{L[A]} ({\delta}) >
	{\aleph_0}$, then $\cf({\delta} ) > {\aleph_0}$. 
\ite 2 $\forces_P $``For all ${\delta} < \chi:$If  $\cf^{L[A]} ({\delta}) >
	{\aleph_0}$, then $\cf({\delta} ) = \aleph_1 $''. 
\ite 3 If $\M' \models \cf^{L[A']}(\mu) > {\aleph_0}$, then $\M'
	\models \cf(\mu ) = \aleph_1 $. 
\ite 4 If $\M' \models ``\cf(\mu) = \aleph_1$'', then $\cf(\mu) = \aleph_1$. 
\ENDent

\proof: 
(1): By the 
choice of $A$.  (\coveringlemma(3)).  
\\
(2):  Use (1) and the fact that $P$ is proper, hence does not cover
	old uncountable sets by new countable sets. 
\\
(3): Use (2) and \absolute. 
\\
(4):  If $\M' \models `` \cf(\mu) = \aleph_1$'', then the function $c'(\mu,
	\cdot)$ is increasing and cofinal in $ \mu$.  (Recall that
	$\omega_1^V  = \omega_1 ^{\M' } $) \QED\thistheorem

\proc Conclusion:  Since $\cf(\mu)={\aleph_0}$, we get from (3) and (4): 
$\M' \models ``L[A'] \models `\cf(\mu) = {\aleph_0}$'\,''.

 Let $\M'
 \models ``\nu_1$ is the $L[A']$-successor of $\mu$.''  We will
 show that $\nu_1=\nu$.   This suffices, because $M'$ is transitive.

So assume that $\nu_1 < \nu$.   Working in $L[A']$ we have
$\bigl|[\mu]^{{\aleph_0}}\bigr| = \nu$ and $|L_{\nu_1}[A']| <\nu$, so we
can find a $y\in
 [\mu]^{\aleph_0} $, $y\in L_{{\gamma} }[A'] \setminus L_{\nu_1}[A']$ for
 some ${\gamma} < \nu$.  Working in $V$, let 
$ L_{\gamma}[A'] = \bigcup_{i < \omega_1} X_i$, where $\< X_i: i <
 \omega _1> $ is a continuous increasing chain of elementary countable
 submodels of $L_{\gamma}[A']$, with $y, A' \in X_0$. In $\M_1' =
 (M_1',\in)$  we can 
 find a continuous increasing sequence $\< Y_i: i < \omega_1 > $ of
 countable elementary submodels of $L_\mu[A']$ 
with $\bigcup_{i< \omega_1} Y_i = L_\mu[A']$ and $A'\in Y_0$.  
 We can find an $i$ such
 that $X_i\cut L_\mu[A']= Y_i$. \\
Let $j: (X_i, \in, A, Y_i) \to (
  L_{\hat{\gamma}}[\hat A],\in , \hat A,  L_{\hat \mu}[\hat A])$ be
  the collapsing 
 isomorphism.  

Now note that $Y_i = X_i \cut L_\mu[A']$ is a transitive subset of
 $X_i$, so $j\on  Y_i $ is exactly the Mostowski collapse of
$(Y_i,\in )$, so $j\on Y_i \in M_1'$ and $\hat A \in M'_1$. Hence
also $j(y) \in 
L_{\hat{\gamma}}[\hat A] \subseteq M_1'$, so we can compute
$$ y = \{ \alpha: (j\on Y_i )(\alpha ) \in j(y) \} $$
in $\M'_1$.  Hence $y\in M'_1 $.   But $\M' \models
``[\mu] ^{{\aleph_0}} \cut M_1' = 
[\mu] ^{{\aleph_0}} \cut M_0' = 
[\mu] ^{{\aleph_0}} \cut L[A'] $'' (the first equality holds because
 $Q_0$ is a $\sigma$-closed forcing notion, the second because of our
 assumption \coveringlemma(3)) 

Hence $\M' \models y \in L[A']$, so $\M' \models y \in L_{\nu_1}[A']$,
 a contradiction to our choice of $\nu$. \qquad\qquad\qquad\null
\QED{\nuisthere{\rm\  Case\ 1}}
\medskip

\noindent{\bf Case 2:} $\cf(\mu)=\aleph_1$.   Let ${\gamma} < \nu$.  We
have to show that ${\gamma}\in M'$.  Since $L[A'] \models
|{\gamma}|=\mu$,  we can in $L[A']$ find an increasing sequence 
$\<A_{\xi}:  {\xi}<\mu>$, ${\gamma} = \bigcup_{{\xi}<\mu }A_{\xi}$,
where each $A_{\xi}$ has (in $L[A']$) 
cardinality $< \mu$.  Let $\alpha_{\xi}$ be the order type of $A_{\xi}$,
then the inclusion map from $A_{\xi}$ into $A_{\zeta} $ 
naturally induces an order
preserving function $f_{{\xi} {\zeta}}: \alpha_{\xi} \to
\alpha_{\zeta}$.    Let  
$B = \< \alpha_{\xi}, f_{{\xi} {\zeta}}: {\xi} \le {\zeta} < \mu>$,
and write $B \on {\beta}$ 
for $ \< \alpha_{\xi}, f_{{\xi} {\zeta}}: {\xi} \le {\zeta} < {\beta}>$. 
Clearly the direct limit of this system is a well-ordered set of order
 type ${\gamma}$.  

So $B$ is in $L[A']$, but we can moreover show that each initial
segment $B\on {\beta} $ 
 is already in $L_\mu[A']$.  This follows from the
fact that each such initial segment can be canonically coded by a
bounded subset of $\mu$.

Since $L_\mu[A'] \subseteq L_\barchi[A'] \subseteq M_1'$, we
know that $B\on {\beta}$ is in $M_1'$ for all ${\beta} < \mu$. 
 In $M_1'$ let
$\< {\xi}_i: i < \omega_1>$ be an increasing cofinal subsequence of
$\mu$.  Let ${\beta}_i = \alpha_{{\xi}_i}$, $h_{i j}= f_{{\xi}_i,
{\xi}_j}$.    Note that the direct limit of the system $\< {\beta}_i,
 h_{i j}; i \le j < \omega_1> $ is still a well-ordered set of order
 type ${\gamma}$.

 So  for each ${\delta} < \omega_1$ we know that the
sequence $b_{\delta} := 
\< {\beta}_i, f_{i j}: i \le j \le {\delta}>$ is in $M_1'$, and
$\M_1' \models b_{\delta}  \in T'$.

Now we can (in $V$) find an uncountable set $C \subseteq \omega_1$ and
 a natural number $n$ such that forall ${\delta} \in C$ we have
 $g'(b_{\delta}) = n$.  Now recall the characteristic property of $g$
 (see \whatisPT) and hence of $g'$ (by \WFandTheory): for each
 ${\delta}_1 < {\delta}_2$ in $C$ we have  a unique branch $B'_\alpha =
 \{b'(\alpha, {\beta}): {\beta} < \omega _1\}$ with $\{b_{\delta_1},
 b_{\delta _2}\} \subseteq B_\alpha$.  A priori this $\alpha$ depends
 on ${\delta}_1$ and ${\delta}_2$, but since $B_\alpha\cut
 B_{{\beta}}=\emptyset $ for $\alpha \not={\beta}$ we must have the
 same $\alpha$ for all ${\delta}\in C$.

So the sequence $\< b_{\delta} : {\delta} \in C> $ is cofinal on some
 branch $B_\alpha'$ which is in $\M'$.  So we get that ${\gamma}$, the
 order type of the limit of this system, is also in $M'$. 
\hbox{\QED{\nuisthere{\rm\ Case\ 2}}%
%		\QED\isacardinal%
	\QED\equicon%
\QED{{\rm[GoSh\ 507]}}%
}

\def\paper#1,{{\it #1}, } 
\def\inbook#1,{in:  #1, }  \def\nextreference#1\par{\bigskip \noindent#1}
\def\journall#1:#2(#3)#4--#5.{#1, {\bf #2} (#3), pp.#4--#5.}
\def\journalp#1:#2(#3)#4.{#1, vol~#2 (#3), p.~#4.}
\def\book#1,{{\it #1},} \def\vol#1 {Vol.~#1}
\def\pages#1--#2{pp.#1--#2}
\def\by#1:{#1, }

\def\AnML{Annals of Mathematical Logic}

\def\AnPAL{Annals of pure and applied logic}

\bigskip

\beginsection
References.

\def\key#1[#2:#3]{\medskip\hangindent=1cm\hangafter1
{\bf [#2]}\write\mgfile{\def\string#1{#2}}%
}

\key\PFApplications[Ba1:PF] J.~Baumgartner, 
	\paper Applications of the Proper forcing axiom, \inbook
	Handbook of set-theoretic topology, 915--959.

 \key\IteratedForcing[Ba2:IF] \by J.~Baumgartner: \paper{Iterated forcing},
        \inbook{Surveys in set theory {\rm
        (A.~R.~D.~Mathias, editor)}}, London Mathematical
        Society Lecture Note Series, No.~8, Cambridge
        University Press, Cambridge, 1983.

\key\fuchino[Fu:PE] Sakae Fuchino, \paper On potential embedding and
	versions of Martin's Axiom, Notre Dame Journal of Formal
	Logic, vol 33, 1992. 

\key \Mitchell[Mi:AT] Bill Mitchell, \paper
Aronszajn Trees and the independence of the transfer property, 
\journall\AnML:5(1971)21--46.

\key\EhrenfeuchtConjecture [Sh 56:EC] Saharon Shelah,  \paper
	Refuting Ehrenfeucht Conjecture on rigid models,   
Proc. of the Symp. in memory of A. Robinson, Yale, 1975, A special 
volume 
in the {\it Israel J. of Math.}, 25 (1976) 273-286.

\key\LXXIII[Sh~73:M2] Saharon Shelah, 
\paper  Models  with  second  order  properties
	II. On   trees  with no undefinable branches, 
	\journall\AnPAL:14(1978)73--87.

%Shelah Proper Forcing, Lecture Notes
 \key\ProperForcing[Sh~b:PF] \by S.~Shelah: \book Proper Forcing, Lecture Notes
        in Mathematics \vol942 , Springer Verlag.    

%Shelah Proper and Improper
 \key\ProperImproper[Sh~f:PI] \by S.~Shelah: \book Proper and Improper
         Forcing, to appear.

\key\todo[To:PF] \by S.~Todor\v cevi\v c: \paper A Note on The Proper
Forcing Axiom, \inbook Axiomatic set theory, J.~Baumgartner,
D.A.~Martin, (eds.), Contemporary Mathematics {\bf 31} (1984).

\bigskip

\obeylines

\baselineskip\normalbaselineskip\parskip0cm

Martin Goldstern               
\smallskip
Institut f\"ur Algebra und Diskrete Mathematik,
%Abteilung fuer Theoretische Informatik        
TU Wien                  
Wiedner Hauptstra{\ss}e 8--10/118.2                 
A-1040 Wien, Austria
{\tt goldstrn@email.tuwien.ac.at}
\smallskip\smallskip
Saharon Shelah
\smallskip
Department of Mathematics
Hebrew University of Jerusalem
91904 Jerusalem, Givat Ram 
Israel
{\tt shelah@math.huji.ac.il}

\bye